\documentclass[a4paper,11pt]{article}

\usepackage[ruled]{algorithm2e}
\usepackage[normalem] {ulem}
\usepackage{color}
\usepackage{amsmath}
\usepackage{graphicx}
\usepackage{amssymb}
\usepackage{subfigure}
\usepackage{a4wide}
\usepackage{authblk}
\usepackage{url}

\newcommand{\unimoreaffil}{\small Department of Physics, Informatics and Mathematics, University of Modena and Reggio Emilia}
\newcommand{\bicoccaaffil}{Department of Informatics, Systems and Communication, University Milano-Bicocca}
\newcommand{\uniboaffil}{Department of Computer Science and Engineering (DISI), University of Bologna}
\newcommand{\ecltaffil}{European Centre for Living Technology, Venezia, Italy}
\newcommand{\daisaffil}{Department of Environmental Sciences (DAIS), University Ca' Foscari}

\newcommand{\mytitle}{Dynamical regimes in non-ergodic random Boolean networks}

\title{\mytitle}

\author[1,2]{Marco~Villani\thanks{Corresponding author (\texttt{marco.villani@unimore.it})\\ The software used for the experiments and the analysis of results presented in this paper is available upon request from the corresponding author.
}}
\author[1]{Davide~Campioli}
\author[4]{Chiara~Damiani}
\author[5]{Andrea~Roli}
\author[2,3]{Alessandro~Filisetti}
\author[1,2]{Roberto~Serra}

\affil[1]{\unimoreaffil}
\affil[2]{\ecltaffil}
\affil[3]{\daisaffil}
\affil[4]{\bicoccaaffil}
\affil[5]{\uniboaffil}

\date{}

\begin{document}

\maketitle

\begin{abstract}
Random boolean networks are a model of genetic regulatory networks that has proven able to describe experimental data in biology. They not only reproduce important phenomena in cell dynamics, but they are also extremely interesting from a theoretical viewpoint, since it is possible to tune their asymptotic behaviour from order to disorder. The usual approach characterizes network families as a whole, either by means of static or dynamic measures. We show here that a more detailed study, based on the properties of system's attractors, can provide information that makes it possible to predict with higher precision important properties, such as system's response to gene knock-out. A new set of principled measures is introduced, that explains some puzzling behaviours of these networks.
These results are not limited to random Boolean network models, but they are general and hold for any discrete model exhibiting similar dynamical characteristics.
\end{abstract}

\section{Introduction}

Biological systems typically involve a great number of interacting units exhibiting a high degree of self-organisation that ensures their continuous functioning and allows them to adaptively respond to environmental changes. Understanding the roots of such nontrivial properties represents a major challenge in theoretical biology and complex systems science.
Among these systems, a prominent role is that of genetic regulatory networks, which in multicellular beings rule the
expression (i.e., the transcription) of the various genes.
A major question related to GRNs behaviour is the coexistence, in multicellular organisms, of several different cell types with a single common genome. It has been proposed~\cite{Kauffman69,Kauffmanhome} that this could be explained by the existence of several possible dynamical behaviours in the same dynamical system. 
In particular, it has been suggested that different cell types correspond to different asymptotic dynamical patterns of activity (attractors) of the same system. According to this view, the genome dictates the units and their basic mutual
interactions, and the attractors determine the cell type. Note that this picture can be expanded to accommodate also epigenetic effects.

A prominent model of genetic regulation is that of Random Boolean networks (RBNs), which were introduced by Kauffman~\cite{Kauffman69,Kauffman71} and became one of the major models of complex systems due to their interesting dynamical behaviour~\cite{Origins,Kauffmanhome,Kadanoff,Bastolla_modular,Bastolla_relevant,aldana_scale_free}. The interest later faded, but in recent years it has been renewed by important theoretical advances~\cite{scaling_socolar,sensitivities,Drossel_number,drossel_review} and also, as far as the application to genetics is concerned, by the availability of genome-wide expression data which can be properly described by RBNs~\cite{serra2004,shmuvelichaldana,serra2007,serra2009genetic,ramo,bornholdt-interface2008}. Further reasons of this renewed interest are related to the possibility to describe complex phenomena, such as cell differentiation~\cite{Villani2010a,Serra2010,Villani2011}, and other biological systems such as whole organisms or tissues~\cite{ingrami,eccs07,acri2008,iet,cheng-interface2008}.

RBNs support  different dynamical regimes and are able to describe the transition from ordered to disordered phases by changing the values of some of their parameters.
It is well known that the {\em bias} (i.e., the {\em internal homogeneity}) and the average number of input gene variables (i.e., the network connectivity) can modulate the order/disorder transition, with higher homogeneity and lower connectivity leading toward more ordered behaviour. In addition, also the choice of the functions associated to each node plays a very important role; indeed, RBN dynamical behaviours are deeply influenced by the various classes of the Boolean functions chosen (see Section~\ref{RBNs} for more details).
It has become commonplace to characterize RBNs as either {\em ordered}, {\em critical} or {\em disordered} depending upon their asymptotic dynamics. Typically, the attractors of ordered networks are stable, in the sense that they can be reached by nearby initial conditions, while disordered networks display a behaviour that is the discrete analogue of the butterfly effect, i.e., a small perturbation of the initial condition usually leads to a different attractor. Because of this latter property, disordered networks are also called {\em chaotic}, with slight abuse of term.
Critical networks display intermediate behaviours and are often considered the ones which can best capture important properties in biological systems.
However, due to their large degree of randomness, the dynamical characterization of these systems is a subtle issue. Commonly, they are characterized on the basis of their structural parameters, thereby providing what is actually a characterization of a whole family (or ensemble) of such networks. But single network realizations may behave in a way very different from the one that is typical of the family. Moreover, even single attractors of the same network can behave in a way that is different from the others. Obviously, the dynamical behaviours of the attractors depend on system's structure, but a common origin does not imply that these dynamical behaviours have to be to similar. 

In this paper, we introduce a detailed way to describe the behaviour of RBNs and we will clarify the relationship between our results and those obtained by the usual ensemble analysis.  We will show that this new set of measures is of great help in clarifying and explaining situations that cannot be properly dealt with the usual methods.\\

The paper is organized as follows: Section~\ref{RBNs} briefly presents the RBN framework; Section~\ref{measure} presents the more common static and dynamic measures of the dynamical RBN regimes, which are deeply discussed in Section~\ref{sec:static-and-dynamics-measures}. Section~\ref{critical} presents our detailed vision of the behaviours observed in the so-called {\em critical} RBNs. In Section~\ref{Consequences} we discuss a particularly significant case in which the new measures allow a new interpretation of experimental biological tests, and Section~\ref{Conclusions} summarises and concludes the work.

\section{Random Boolean network models}
\label{RBNs}
In this section we succinctly review RBNs, emphasizing the most relevant properties and results for this work. The interested reader is referred to existing reviews for a more thorough discussion of
RBNs~\cite{Kauffman69,Origins,Kauffmanhome,Kadanoff,drossel_review}.

A Boolean network (BN) is a discrete-state and discrete-time dynamical system whose structure is defined by a directed graph of $N$ nodes, each associated to a Boolean variable $x_i$, $i = 1, \ldots, N$, and a Boolean function $F_i(x_{i_1}, \ldots, x_{i_{k_{in,i}}}$), where $k_{in,i}$ is the number of inputs of node $i$. The arguments of the Boolean function $F_i$ are the values of the nodes whose outgoing arcs are connected to node $i$. The state of the system at time $t$, $t \in \mathbb{N}$, is defined by the array of the $N$ Boolean variable values at time $t$: $s(t) \equiv (x_1(t), \ldots, x_N(t))$.
The most studied BN models are characterized by having a \emph{synchronous} dynamics---i.e., nodes update their states at the same instant---and \emph{deterministic} functions. However, many variants exist, including asynchronous and probabilistic update rules~\cite{pbn}.

A special category of BNs that has received particular attention is that of RBNs, which can capture relevant phenomena in genetic and cellular mechanisms and complex systems in general.
RBNs are usually generated by choosing at random $k_{in}$ inputs per node and by defining the Boolean functions by assigning to each entry of the truth tables a 1 with probability $p$ and a 0 with probability $1 - p$. The parameter $p$ is referred to as the \textit{bias}. Depending on the values of $k_{in}$ and $p$ the dynamics of RBNs is called either \textit{ordered} or \textit{chaotic}.
In the first case, the majority of nodes in the attractor is frozen and any moderate-size perturbation is rapidly dampened and the network returns to its original attractor.
Conversely, in chaotic dynamics, attractor cycles are very long and the system is extremely sensitive to small perturbations: slightly different initial states lead to divergent trajectories in the state space.
RBNs temporal evolution undergo a second order phase transition between order and chaos, governed by the following relation between $k_{in}$ and $p$: 
\begin{equation}
\label{criticalityEquation}
k_{in}^c = [2 p_c (1 - p_c)]^{-1}
\end{equation}
\noindent
where the superscript $c$ denotes the critical values~\cite{poumeau}.

\section{The measure of RBN dynamical regimes}
\label{measure}

In this section we first introduce and discuss the usual ways for characterizing the dynamical regimes of RBNs. Then, we revise the notion of {\em critical network} showing that the dynamic behaviour of a BN can be dramatically different across its attractors.

\subsection{The spread of perturbations}
As previously observed, ordered and disordered dynamical regimes are usually described as systems' behaviours supporting respectively short attractors with fairly regular {\em basins of attraction},\footnote{An attractor basin is the set of states whose evolution lead to the attractor, its size (or dimension) being the cardinality of the set.} and long attractors with sensitive dependence upon initial conditions.\footnote{We remind that in disordered systems trajectories starting from nearby points typically lead to different attractors.}
There are two different ways to identify ordered and disordered regions: \emph{(i)} a static one, based upon the knowledge of the nodes average connectivity and the bias of the Boolean functions and \emph{(ii)} a dynamical one based upon the study of the spreading of perturbations through the system.

In ordered networks a change at one node (i.e., a bit flip) propagates in one step on average to less than one other node: starting from a random initial condition an ordered system rapidly reaches a stable condition in which the majority of nodes are frozen; if this asymptotic behaviour is perturbed, there is a very high probability of coming back to it. On the contrary, in disordered networks a perturbation at one node propagates in one time step on average to more than one other node: very close initial conditions could rapidly diverge toward different attractors, and attractors typically have long periods, with large portions of oscillating nodes; if perturbed, the system has therefore high probability of changing its asymptotic behaviour. Critical systems are at the boundary between these two dynamical regimes: a change at one node propagates in one time step on an average to exactly one other node. This situation corresponds to the percolation of frozen nodes through the network, and therefore to the formation of wide but isolated oscillating zones.

\subsection{Static vs. dynamic estimates}
The main static methods to measure the RBN dynamical regimes implicitly presume {\em ergodicity}, that is, inputs arise with the same probability during evolution, and time average over the states visited by the network yields the same results as averages over the whole phase space. 
A widely used measure is the {\em average sensitivity} of a Boolean function, proposed by Shmulevich and Kauffman~\cite{sensitivities}. The {\em sensitivity} of a function $F_i$ measures how sensitive the output of the function is to changes of its inputs. Let us consider an input configuration $x$ and all its 1-Hamming neighbours, i.e., the input configurations that differ from $x$ in exactly one value; the sensitivity $s^{F_i}(x)$ is defined as the number of 1-Hamming neighbours of $x$ on which the $F_i$ function values are different than on $x$~\cite{sensitivities}.
The {\em average sensitivity} of the function $F_i$, $I(F_i)$ is the expected value of sensitivity $s^{Fi}(x)$ with respect to the distribution of $x$.
Related to the average sensitivity is the notion of {\em influence} of a variable: the {\em influence} of the $j$-th input variable of a function $F_i, I_j(F_i)$, is the probability that the function $F_i$ changes its value when the value of the $j$-th variable is changed (a concept linked to Boolean derivatives and Lyapunov exponents of RBNs~\cite{Bagnoli1992,luque}).
It can be proven that $I(F_i)$ is $k_{in}$ times the average probability that the output of a node changes when one of its inputs changes, in formulas:
\begin{equation}
\label{sensitivityEquation}
I(F_i)=(1-q_i)k_{in,i}
\end{equation}
where $q_i$ is the probability that the output of the $i-th$ node does not change when one of its $k_{in,i}$ inputs changes~\cite{sensitivities}.

At this point Shmulevich and Kauffman introduce the ergodic hypothesis: under uniform input distribution the average sensitivity of function $F_i$ is simply equal to the sum of the influences of all its input variables (a number that ranges from 0 to $k_{in}$). The average sensitivity of a BN is the weighted average of the sensitivities of all its functions. If this final average sensitivity is lower, equal to or higher than 1 the network typically dampens, maintains or amplifies the perturbations and is therefore ordered, critical or disordered, respectively. It should  be noted that, when the Boolean functions are drawn according to a Bernoulli distribution with a specific bias $p$, the average sensitivity coincides with the well-known critical transition curve $k_{in}(2 p (1-p))$ (see~\cite{sensitivities}). However, when network functions are generated according to probability distributions that favour some variables relative to others, or when functions are chosen randomly from certain classes of functions (e.g., canalizing), the average sensitivity well captures the dynamical behaviour of the network.
We also remind that the ergodicity assumption does not hold for the dynamics of arbitrary RBNs~\cite{canalizing,Szejka2008}.

The alternative to static measures is that of explicitly exploring the dynamical behaviour of the system. An interesting and well-known method is the so-called Derrida procedure~\cite{poumeau,derrida86} which  exploits the spreading of perturbations through the network. This procedure involves two parallel runs of the same system, with initial states that differ for only a small fraction of the units. This difference is usually measured by means of the Hamming distance. Let $h(t)$ be the Hamming distance between the states at time $t$ in the two parallel trajectories. If after a transient the two runs are likely to converge to the same state,\footnote{The estimation is performed on many different initial conditions.} i.e., $h(t)\rightarrow0$, then the dynamics of the system is robust with respect to small perturbations---a signature of the ordered regime. Conversely, if the difference is likely to increase, then the dynamics is sensitive to small perturbations and the corresponding regime is disordered. 
The common way to measure the dynamical regime of a RBN by means of the Derrida procedure is that of randomly generating a great number of pairs of initial conditions differing for one or a few units, evolving the network for one step, measuring the Hamming distance of the two resulting states, taking the averages for each perturbation size, and plotting these averages (i.e., $h(1)$) vs. the initial perturbation size (i.e., $h(0)$). 
Let $\lambda$ be the slope at the origin of the curve $h(1)$ vs. $h(0)$.
A system is ordered for $\lambda < 1$, whereas it is disordered for $\lambda > 1$. The case of $\lambda = 1$ identifies the critical regime~\cite{Kadanoff}.
It is worth stressing that the Derrida parameter is an empirical estimation of the average sensitivity mentioned above. Furthermore, under the hypothesis of ergodicity, it holds $\lambda  \sim (1-q)k_{in}$, where $q$ is the average probability that the output of a node does not change when one of its inputs does~\cite{sensitivities}.

\subsection{Attractor sensitivity}
It is important to notice that the usual interpretation of the behaviour of RBNs, e.g., as models of genetic regulatory networks, privileges their attractor states. Indeed, the network will be found in one of these states, after transients have died out. Therefore, measures taken on randomly chosen states do not necessarily provide a correct estimate of the relevant dynamical behaviour of the system, which should rather be evaluated on its attractors.

A more meaningful way to analyse the dynamics is that of applying the classical Derrida procedure only on the states belonging to the attractors~\cite{campioli2011,VillaniSerra-attrPert-2014}. We can therefore define the {\em sensitivity on attractor} $i$ ($SA_i$) as the result of the Derrida procedure performed only on the states belonging to the attractor $i$, and the {\em sensitivity on attractors} ($SA$) as the average of the $SA_i$, each $SA_i$ being weighted with the size of its attraction basin. We will also refer to the usual Derrida procedure, i.e., the one on random initial states, as $DA$.

A similar approach was used in~\cite{drossel_generalized_derrida}, where the authors only considered what has been defined here as the quantity $SA$. As we will see in the following, the same RBN can support several attractors with different dynamical behaviour, therefore the $SA$ value alone does not provide a complete picture and, in some conditions, can even suggest misleading conclusions on the dynamical regime.
In order to appreciate differences and similarities among $DA$, $SA$ and $SA_i$, we will analyse these measures across different statistical ensembles.
In the following we use the word {\em family} for networks with transition functions randomly picked from the same set $\mathcal{F}$. Families will be denoted by the symbols $\mathcal{M}_1,\mathcal{M}_2,\dots,\mathcal{M}_i$.
We assume that each node has the same number of incoming links $k_{in}$, and that the origin of these links is chosen at random with uniform probability among the remaining $N-1$ nodes, prohibiting multiple connections. As previously mentioned, in random networks by varying $\mathcal{F}$ and/or $p$ it is possible to move from disordered to ordered dynamical regimes.

\begin{figure}[t]
\centering
\includegraphics[width=0.8\textwidth]{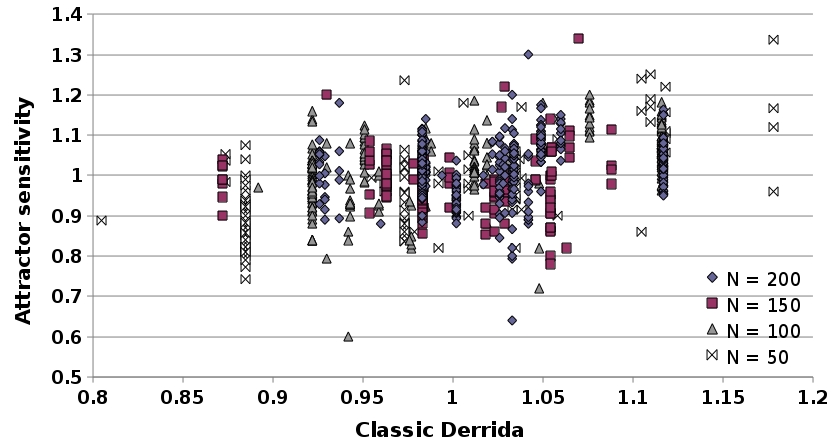}
\caption{Attractor sensitivity of the attractors of 40 networks. The networks are grouped in 4 sets having different number of nodes. As it can be noted, the attractor sensitivities of nets having the same $DA$ (vertical sequences of marks) can have very different values.}
\label{DAvsSAi}
\end{figure}

\begin{figure}[t]
\centering
\includegraphics[width=0.6\textwidth]{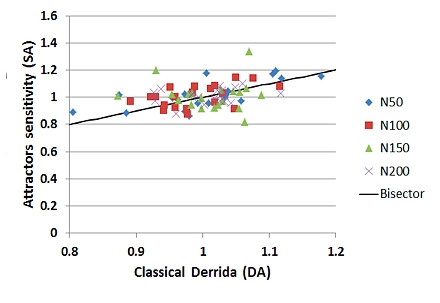}
\caption{The $SA$ values of the same nets of Figure~\ref{DAvsSAi}. $SA$ correlates with $DA$ so that there is a clear---although noisy---proportionality between $SA$ and $DA$. The figure shows also the bisector (the continuous line).}
\label{DAvsSA}
\end{figure}

Let us first consider the usual ensemble, i.e., families composed of networks of different size with random topology, $k_{in}=2$ and Boolean functions uniformly distributed ($p=0.5$), referred to as $\mathcal{M}_1$ family in the following.
As one can observe in Figure~\ref{DAvsSAi}, the attractors of the same network can have values of $SA_i$ that significantly differ from one another and that in turn may be very different from the corresponding $DA$ value of the network. These differences can be observed in several networks across different families characterized by different number of nodes. However, one can also see from Figure~\ref{DAvsSA}a that $SA$ (i.e. the weighted average of the attractors' sensitivities) is strongly correlated with $DA$ (note that the relation is well approximated by a linear function of slope 1).
This result means that the dynamical behaviour averaged across all the attractors of a single RBN is related to the dynamical behaviour of a random sampling of the systems' state space. 

\begin{table}[t]
\centering
\caption{The averages of $DA$ and of $SA$ for the 4 sets of nets in Figure~\ref{DAvsSAi} and~\ref{DAvsSA}.\label{tab:DAvsSA}}{%
\begin{tabular}{|r|r|r|}
\hline
\textbf{N} & \textbf{$\langle DA \rangle$} & \textbf{$\langle SA \rangle$} \\
\hline
50 & 1.00$\pm$0.09 & 1.0$\pm$0.1\\
\hline
100 & 0.99$\pm$0.06 & 1.00$\pm$0.05\\
\hline
150 & 1.00$\pm$0.05 & 1.01$\pm$0.05\\
\hline
200 & 1.00$\pm$0.05 & 1.01$\pm$0.05\\
\hline
\end{tabular}
}
\end{table}

Even more remarkably, from the results in Table~\ref{tab:DAvsSA} we observe that by averaging over networks with the same size, the average values of $DA$ and of $SA$ are almost equal.\footnote{In this case, the averages are around 1, as the ensemble considered is the first historical example of critical systems.} 
A similar situation holds for other families generated with different values of $k_{in}$ and $p$, in all the three dynamical regimes.\marginpar{} 
In all these families, the attractor sensitivity considerably differs even across attractors of the same network. However, $DA \approx SA$ when the values are averaged across all the networks of the family.

It is important to observe that this last property does not hold for all the possible families. Let us now consider two families of RBNs---in the following denoted by $\mathcal{M}_5$ and $\mathcal{M}_6$, respectively---with $k_{in}=2$, in which only a subset of canalizing Boolean functions is allowed. The interest for these particular families comes from a study on random threshold networks (RTNs), in which each node computes the weighted sum of its input (we consider the case of weights in $\{-1,+1\}$). If the sum exceeds a given threshold, the node takes the value 1; 0, otherwise.
It is possible to map these transition rules on the RBN framework, obtaining a set $\mathcal{F}$ of rules for each threshold. In this work, we study the dynamical behaviour of the RBN families $\mathcal{M}_5$ and $\mathcal{M}_6$ , whose activation patterns correspond to RTNs having respectively thresholds $h=+0.5$ and $h=-0.5$.
Incidentally, the two sets of Boolean functions identified in this way are complementary to each other. Despite this similarity, the asymptotic behaviours of these two families are quite different, and their $SA$ value is considerably different from the classical Derrida measure---see Table~\ref{tablem5m6} for the details. 
Indeed, the attractor sensitivity is $0.96$ for the $\mathcal{M}_5$ family and  $0.65$ for the  $\mathcal{M}_6$ family (70 nodes per network). This difference can provide an explanation for the striking different dynamical behaviour in the two families: the average number of the attractors of networks in $\mathcal{M}_5$ is significantly higher than the average number of attractors in $\mathcal{M}_6$; considerable differences can also be found in length of cycles and number of frozen nodes~\cite{campioli2011}. The two families of RBNs show therefore very different behaviours despite the fact that both have the same $DA$. Hence, the Derrida parameter is not able to correctly describe the dynamics in these cases, where there is a remarkable difference between $DA$ and $SA$. 

\begin{table}[t]
\centering
\caption{The Boolean functions allowed in the $\mathcal{M}_5$ and $\mathcal{M}_6$ families ($N=70$, $k_{in}=2$) and their respective measures of dynamical behaviours.\label{tablem5m6}}{%
\resizebox{\textwidth}{!}{%
\begin{tabular}{|cc|cccc|cccc|}
\cline{3-10}    \multicolumn{2}{c|}{  } & \multicolumn{4}{c|}{ $\mathcal{M}_5$}       & \multicolumn{4}{c|}{ $\mathcal{M}_6$} \\
    \hline
   \textbf{$A$} & \textbf{$B$} & \textbf{$A \vee B$} & \textbf{$\neg A \wedge B$} & \textbf{$A \wedge \neg B$} & \textbf{$FALSE$} & \textbf{$\neg(A \vee B)$} & \textbf{$\neg A \vee B$} & \textbf{$A \vee \neg B$} & \textbf{$TRUE$} \\
  
    \hline
    0     & 0     & 0     & 0     & 0     & 0     & 1     & 1     & 1     & 1 \\
    \hline
    0     & 1     & 1     & 1     & 0     & 0     & 0     & 1     & 0     & 1 \\
    \hline
    1     & 0     & 1     & 0     & 1     & 0     & 0     & 0     & 1     & 1 \\
    \hline
    1     & 1     & 1     & 0     & 0     & 0     & 0     & 1     & 1     & 1 \\
    \hline
    \multicolumn{2}{|c|}{\textbf{{average sensitivity\footnote{Value derived from Equation~\ref{sensitivityEquation}\hfil}}}} & \multicolumn{4}{c|}{0.75}     & \multicolumn{4}{c|}{0.75} \\
    \hline
    \multicolumn{2}{|c|}{\textbf{{DA}}} & \multicolumn{4}{c|}{0.74}     & \multicolumn{4}{c|}{0.74} \\
    \hline
    \multicolumn{2}{|c|}{\textbf{{SA}}} & \multicolumn{4}{c|}{0.96}     & \multicolumn{4}{c|}{0.67} \\
    \hline
    \end{tabular}
}}
\end{table}%

\begin{figure} 
\centering
 \includegraphics[width=1\textwidth]{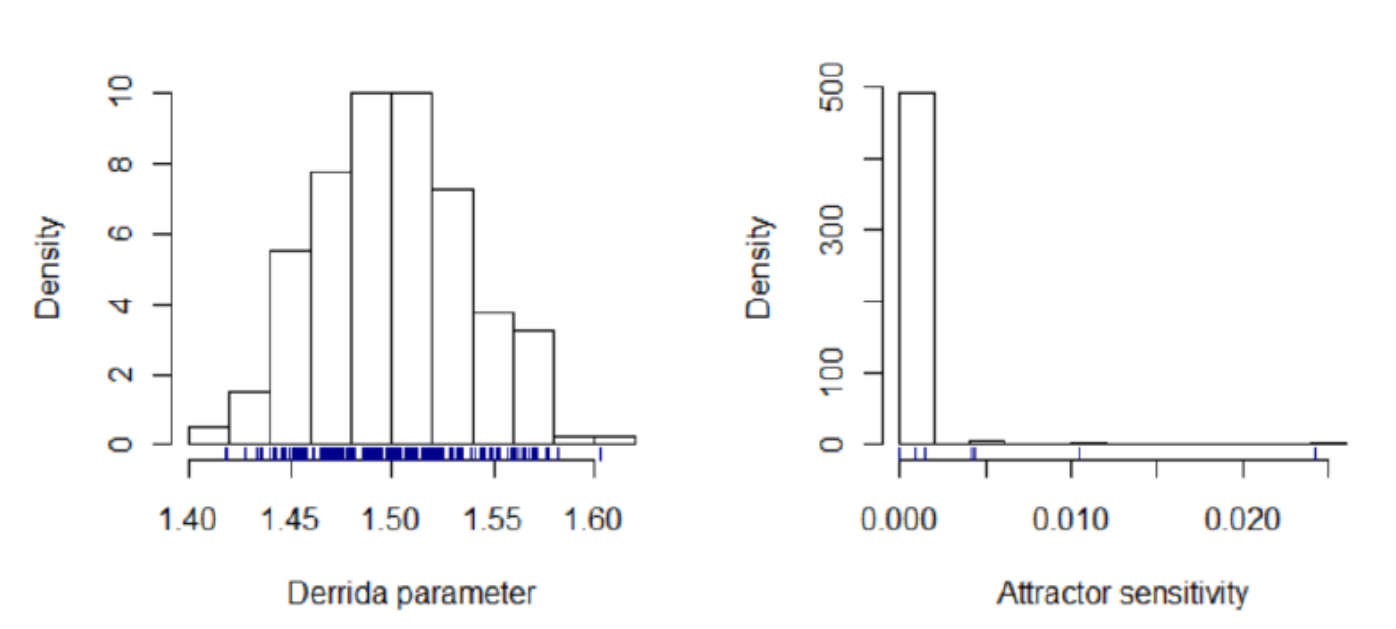}
  \caption{Probability distribution of classical Derrida values ($DA$) and of attractor sensitivity ($SA$) for RBNs of 71 nodes in the family $\mathcal{M}_7$. 200 networks sampled.}
  \label{SAevolute}
\end{figure}

An even more extreme case in which $DA$ and $SA$ are very different concerns the case of BNs evolved to perform a global computation task. We consider the case of the so-called {\em Density Classification Problem}, which consists of classifying a binary string in either of two classes, depending on the ratio between 0s and 1s. The goal is then to find the Boolean functions of a network such that it is driven to a uniform state, consisting of all 1s, if the initial configuration contains more 1s and all 0s otherwise~\cite{packard}. 
It can be shown that the simple {\em majority rule} applied on random topologies outperforms all human or artificially-evolved rules running on ordered lattice, with a performance that increases with the size of the network~\cite{serra2002,Mesot1920}. The majority rule states that the value of a node at time $t+1$ is 0 (resp. 1) if the majority of its neighbours has value 0 (resp. 1) at time $t$.
In this context, we studied a family of BNs ($\mathcal{M}_7$), with $k_{in}=3$ and $N=71$,\footnote{An odd number of nodes avoids the cases with equal quantity of 0s and 1s} evolved to solve this task (see~\cite{Benedettini2013} for details). We are interested here in the analysis of the evolved networks in terms of sensitivity and Derrida parameter. According to the value of $DA$, this family is chaotic (with $DA=1.50$), but the attractor sensitivity $SA$ supports a very different conclusion: we have $SA < 0.001$, indicating that the system is deeply in the ordered region. And indeed the system is very ordered, having just a few very short attractors (typically, only two fixed points) with regular basins of attraction (nearby initial conditions typically evolve towards the same attractor). Figure~\ref{SAevolute} depicts the distributions of $DA$ and $SA$ in 200 realizations, showing the striking difference existing in this case between these parameters.

In summary, we remark that in general $DA$ and $SA$ can considerably differ. Furthermore, even when they coincide, we must be very careful in drawing conclusions on the dynamical regimes, as the individual sensitivities on single attractors can be very dissimilar. In the next section we discuss the reasons for these differences and we outline their implications.

\section{Linking static and dynamical measures}
\label{sec:static-and-dynamics-measures}

The point where static and dynamical measures deviate is that of the assumption of ergodicity made by the static approach. As previously pointed out, under the hypothesis of uniform input distribution, it can be shown that the average sensitivity of a function $F_i$ is equal to the sum of the influences of all its input variables~\cite{sensitivities}. Nevertheless, in all the other cases this computation requires more attention. In particular, the asymptotic states sample only a part of the whole state space, a part that in non-ergodic systems can be quite limited. In order to circumvent this problem one needs a way to correctly weight each input configuration.
In RBNs it is possible to estimate the $SA_i$ of the $i$-th attractor from its time series by computing the average fraction of ``1'' (i.e., its frequency of occurrence $b$), and weight in such a way the influence of the input variables of each function. To obtain the average sensitivity of function $F_i$, i.e., $I(F_i)$, we sum these weighted influence values and then we take the global average over all the Boolean functions, weighting each function by its occurrence in the network (Table~\ref{averageSensitivities}). The global measure $SA$ is computed by averaging over all the attractor sensitivity values, $SA_i$, weighting each $SA_i$ with the basin of attraction size of the corresponding attractor. Table~\ref{BooleanFunctions} and Table~\ref{averageSensitivities} show the values of a particular $SA_i$ for $\mathcal{M}_5$ and $\mathcal{M}_6$ families, whereas Table~\ref{dependenceSensitivities} shows the dependence of $SA$ from the number of nodes $N$.

We observe that $b$, i.e., the asymptotic number of ``1'', depends on the dynamics of the system, which in turn originates from system's features such as connectivity and Boolean functions. Therefore, it should be possible to analytically estimate this value. 
Indeed, one may use the so-called {\em annealed approximation}, which is a sort of mean field approximation---that holds for annealed networks~\cite{poumeau} with an infinite number of nodes---to estimate $b$ and hence the value of $SA$. This approach can give reasonable guesses also for quenched systems~\cite{poumeau,Szejka2008}. 
Nevertheless, it takes averages and therefore ignores the possible peculiar behaviour on individual attractors; moreover, it cannot take into account the effect of the finite number of nodes that may be often non-negligible. An example of this discrepancy can be observed for the family $\mathcal{M}_5$ (see Table~\ref{dependenceSensitivities}).
Eventually, by knowing the asymptotic sensitivity $I(F_j)$ of function $F_j$ and by using Equation~\ref{sensitivityEquation}, it is possible to derive for each node $j$ the probability of output change in case of an input change (or the complementary probability $q_j$ of maintaining its value constant). 

\begin{table}[htbp]
  \centering
\caption{The table shows the probability of occurrence $P(A,B)$ of each input configuration vs. $b$, i.e., the normalized occurrence of value 1. A ``CA'' is reported if a change on the flipped input has as a consequence a change in the corresponding function output.\label{BooleanFunctions}}{
\resizebox{\textwidth}{!}{%
\begin{tabular}{|r|r|r|cccc|cccc|}
    \hline
    \multicolumn{2}{|r|}{\textbf{Input}} & \textbf{$P(A,B)$} & \multicolumn{4}{c|}{\textbf{Family $\mathcal{M}_5$}} & \multicolumn{4}{c|}{\textbf{Family $\mathcal{M}_6$}} \\
    \hline
    \multicolumn{1}{|c}{\textbf{$A*$}} & \multicolumn{1}{c|}{\textbf{$B$}} & \textbf{} & \textbf{$A \vee B$} & \textbf{$A \wedge \neg B$} & \textbf{$\neg A \wedge B$} & \textbf{$FALSE$} & \textbf{$\neg (A \vee B)$} & \textbf{$\neg A \vee B$} & \textbf{$A \vee \neg B$} & \textbf{$TRUE$} \\
    \hline
    \multicolumn{1}{|c}{0} & \multicolumn{1}{c|}{0} & $(1-b)^2$ & CA    & CA    & nc    & nc    & CA    & CA    & nc    & nc \\
    \hline
    \multicolumn{1}{|c}{0} & \multicolumn{1}{c|}{1} & $b(1-b)$ & nc    & nc    & CA    & nc    & nc    & nc    & CA    & nc \\
    \hline
    \multicolumn{1}{|c}{1} & \multicolumn{1}{c|}{0} & $b(1-b)$ & CA    & CA    & nc    & nc    & CA    & CA    & nc    & nc \\
    \hline
    \multicolumn{1}{|c}{1} & \multicolumn{1}{c|}{1} & $b^2$    & nc    & nc    & CA    & nc    & nc    & nc    & CA    & nc \\
    \hline
    \hline
\multicolumn{1}{|c}{\textbf{$A*$}} & \multicolumn{1}{c|}{\textbf{$B$}} & \textbf{} & \textbf{$A \vee B$} & \textbf{$A \wedge \neg B$} & \textbf{$\neg A \wedge B$} & \textbf{$FALSE$} & \textbf{$\neg (A \vee B)$} & \textbf{$\neg A \vee B$} & \textbf{$A \vee \neg B$} & \textbf{$TRUE$} \\
    \hline
    \multicolumn{1}{|c}{0} & \multicolumn{1}{c|}{0} & $(1-b)^2$ & CA    & nc    & CA    & nc    & CA    & CA    & CA    & nc \\
    \hline
    \multicolumn{1}{|c}{0} & \multicolumn{1}{c|}{1} & $b(1- b)$ & CA    & nc    & CA    & nc    & CA    & nc    & CA    & nc \\
    \hline
    \multicolumn{1}{|c}{1} & \multicolumn{1}{c|}{0} & $b(1- b)$ & nc    & CA    & nc    & nc    & nc    & CA    & nc    & nc \\
    \hline
    \multicolumn{1}{|c}{1} & \multicolumn{1}{c|}{1} & $b^2$    & nc    & CA    & nc    & nc    & nc    & CA    & nc    & nc \\
    \hline
    \end{tabular}%
}}
\end{table}%

\begin{table}[htbp]
 \centering
\caption{The table shows the influence of the input variables, computed by taking into account the effective occurrence probability of the possible input configurations, and the resultant function sensitivity. The sensitivity of the considered attractor, $SA_i$, is estimated (by construction, in each family the four Boolean functions have the same occurrence probability).\label{averageSensitivities}}{
\resizebox{\textwidth}{!}{%
   \begin{tabular}{|r|r|r|r|r|r|r|r|r|}
    \hline
          & \multicolumn{4}{r|}{\textbf{Family $\mathcal{M}_5$ ($b=0.08$)}} & \multicolumn{4}{r|}{\textbf{Family $\mathcal{M}_6$ ($b=0.67$)}} \\
    \hline
          & \textbf{$A \vee B$} & \textbf{$A \wedge \neg B$} & \textbf{$\neg A \wedge B$} & \textbf{$FALSE$} & \textbf{$\neg (A \vee B)$} & \textbf{$\neg A \vee B$} & \textbf{$A \vee \neg B$} & \textbf{$TRUE$}  \\
    \hline
    \textbf{Influence of A} & 0.92  & 0.08  & 0.92  & 0     & 0.33  & 0.22  & 0.67  & 0  \\
    \hline
    \textbf{Influence of B} & 0.92  & 0.92  & 0.08  & 0     & 0.33  & 0.78  & 0.33  & 0  \\
    \hline
    \textbf{Function sensitivity:} & 1.84  & 1     & 1     & 0     & 0.66  & 1     & 1     & 0  \\
    \hline
    \textbf{$SA_i$:} & 0.96  &       &       &       & 0.67  &       &       &   \\
    \hline
    \end{tabular}
}}
\end{table}%

\begin{table}[htbp]
  \centering
\caption{For each family the table shows the theoretical estimate of $SA$ (computed by means of the procedure explained in the text) and the corresponding experimental value (estimated by effectively performing the Derrida procedure on random initial conditions and only on the attractor states respectively for $DA$ and $SA$ measures), vs. different number of nodes.\label{dependenceSensitivities}}{
    \begin{tabular}{|r|r|r|r|r|r|r|r|}
    \hline
    \textbf{Measure} & \textbf{N} & \textit{\textbf{b}} & \multicolumn{2}{r|}{\textbf{Family  $\mathcal{M}_5$ }} & \textit{\textbf{b}} & \multicolumn{2}{r|}{\textbf{Family  $\mathcal{M}_6$ }}  \\
    \hline
          &       & \textit{} & \textbf{Theo.} & \textbf{Exper.} & \textbf{} & \textbf{Theo.} & \textbf{Exper.}  \\
    \hline
    \textbf{DA} & 70    & 0.5   & 0.75  & 0.74  & 0.5   & 0.75  & 0.74  \\
    \hline
    \textbf{SA} & 70    & 0.08  & 0.96  & 0.93  & 0.67  & 0.66  & 0.64  \\
    \hline
    \textbf{SA} & 700   & 0.05  & 0.97  & 0.95  & 0.66  & 0.67  & 0.67  \\
    \hline
    \textbf{SA} & $\infty$     & 0     & 1     & ---   & 0.67  & 0.67  & ---  \\
    \hline
    \end{tabular}%
}
\end{table}%

\begin{figure} 
\centering
 \includegraphics[width=1\textwidth]{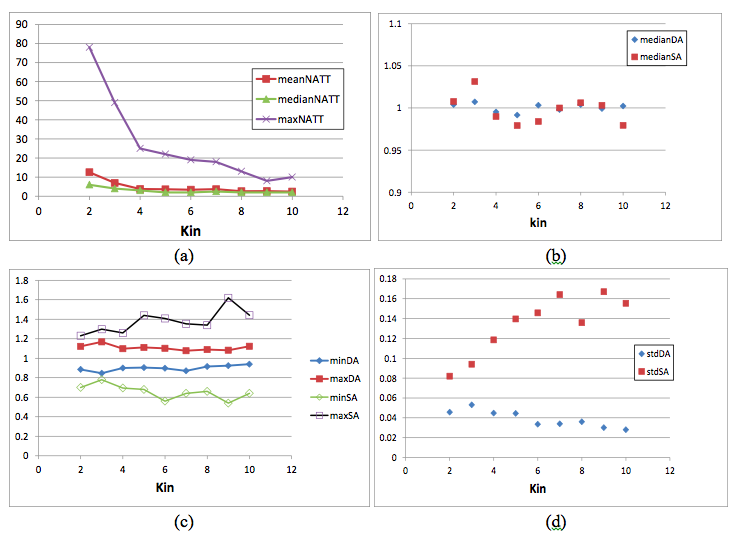}
  \caption{Measures on ``critical'' families (their bias $p$ deriving from Equation~\ref{criticalityEquation} with $N=200$ and $k_{in}$ variable from 2 to 10, 50 networks for each family. {\em (a)} Number of attractors in each RBN (average, median, maximum); {\em (b)} median values of $DA$ and $SA$; {\em (c)} minimum and maximum value of $DA$ and $SA$; {\em (d)} the dispersion (standard deviation) of the same distribution in {\em (b)} and {\em (c)}.}
  \label{reticritiche}
\end{figure}

\section{``Critical'' networks can have very different dynamical behaviours}
\label{critical}
The aforementioned considerations provide us with instruments to appreciate in finer details the differences among the networks usually classified ``at the edge of chaos''~\cite{Origins}: the values of their static parameters are all ``critical'', but the dynamical behaviours of their attractors could significantly deviate from critical regimes. This phenomenon finds striking evidence in the results of a series of experiment we made, in which we created 50 networks for $k_{in} \in \{2,3,\ldots,10\}$ and $N=200$, with Boolean functions generated by using the bias $p = p_c$, where $p_c$ is the critical value derived from Equation~\ref{criticalityEquation}.
In Figure~\ref{reticritiche}b, we can see that the statistical ensembles show a critical behaviour, as the measures $DA$ and $SA$ assume values very close to 1. Nevertheless, the number of attractors in each RBN class dramatically decreases as $k_{in}$ increases (Figure~\ref{reticritiche}a). Furthermore, the minimum and maximum values of $DA$ remain constant as $k_{in}$ increases, whereas the corresponding minimum and maximum values of $SA$ increase their (already remarkable) distance from 1.0 (see Figure~\ref{reticritiche}c); this fact can also be observed looking at the dispersion of their respective distributions (the standard deviations of $DA$ and $SA$ distributions shown in~Figure~\ref{reticritiche}d). Obviously, the single attractor sensitivities $SA_i$ are still more extreme~\marginpar{}. Therefore, although the statistical ensembles maintain a critical behaviour across the values of $k_{in}$, each single RBN realization has higher and higher probability to deviate from the critical behaviour; furthermore, the magnitude of this deviation grows as $k_{in}$ increases.
Interestingly, the same deviation increase can be observed also in RBNs whose bias $p$ differs from the ``edge of chaos'' values; for example, we observed ordered RBN ensembles in which, by increasing $k_{in}$, one observes that there are progressively more and more RBN realizations with differing $DA$ and $SA$, although their ensemble averages remain close to each other. 
On the basis of these results, we can claim that by using these measures it is possible to introduce a finer classification within the ``classical'' critical networks that allows a clearer understanding of their asymptotic dynamical behaviour. 
In addition, it is worth to emphasize that the discrepancy between this finer analysis and the usual one increases as $k_{in}$ increases.

\section{Consequences for biological behaviours}
\label{Consequences}
The previous considerations have important consequences on the application of RBNs as model of genetic regulatory networks.  It is very interesting that, on the same network realization, $DA$ and $SA$ can have different values, but even more remarkable is the fact that the same RBN can support attractors showing different attractor sensitivities and therefore different dynamical behaviours, that can hopefully be accessible to experimental measures.
As a leading example, we can use a prominent interpretation of biological measures, performed by using the RBN framework to describe gene knock-out experiments in {\em Saccharomyces cerevisiae}~\cite{serra2004,serra2007,graudenzi_scalefree}.

A gene knock-out is a permanent silencing of a single gene: biologists can experimentally induce such modification and quantify (by using for example cDNA microarrays) the corresponding change of expression of all the other genes (the so-called avalanche)~\cite{hughes}.
It is possible to reproduce {\em in silico} the knock-out process: at a certain time point, starting from a state of an attractor,\footnote{for the reasons discussed above} the value of one of the nodes of a RBN is permanently clamped to the value 0. The evolution of the unperturbed network is compared to that of the perturbed one (after the latter has reached again a stable asymptotic behaviour). A gene is said to be affected (or perturbed) if its behaviour differs in the two networks in at least one time step. The avalanche corresponding to a given knock-out is the set of perturbed genes (including the one which has been knocked-out), and its size is the number of perturbed genes.
We have shown previously~\cite{serra2004,serra2007} that {\em Saccharomyces cerevisiae} shows a dynamical behaviour close to the ``edge of chaos'', on the side of the ordered region and that for random topologies the probability of having an avalanche involving $m$ genes is given by the formula:
\begin{equation}
\label{valanga}
P_m=B_m\lambda^{m-1}e^{-m\lambda}
\end{equation}
where $B_m$ is a coefficient\footnote{See~\cite{serra2007} for details on the calculation of these coefficients.} and $\lambda$ is equal to $(1-q) \langle k_{out} \rangle$, where $\langle k_{out} \rangle$ is the average of the number of outgoing connections per node.

\begin{figure} 
\centering
 \includegraphics[width=0.8\textwidth]{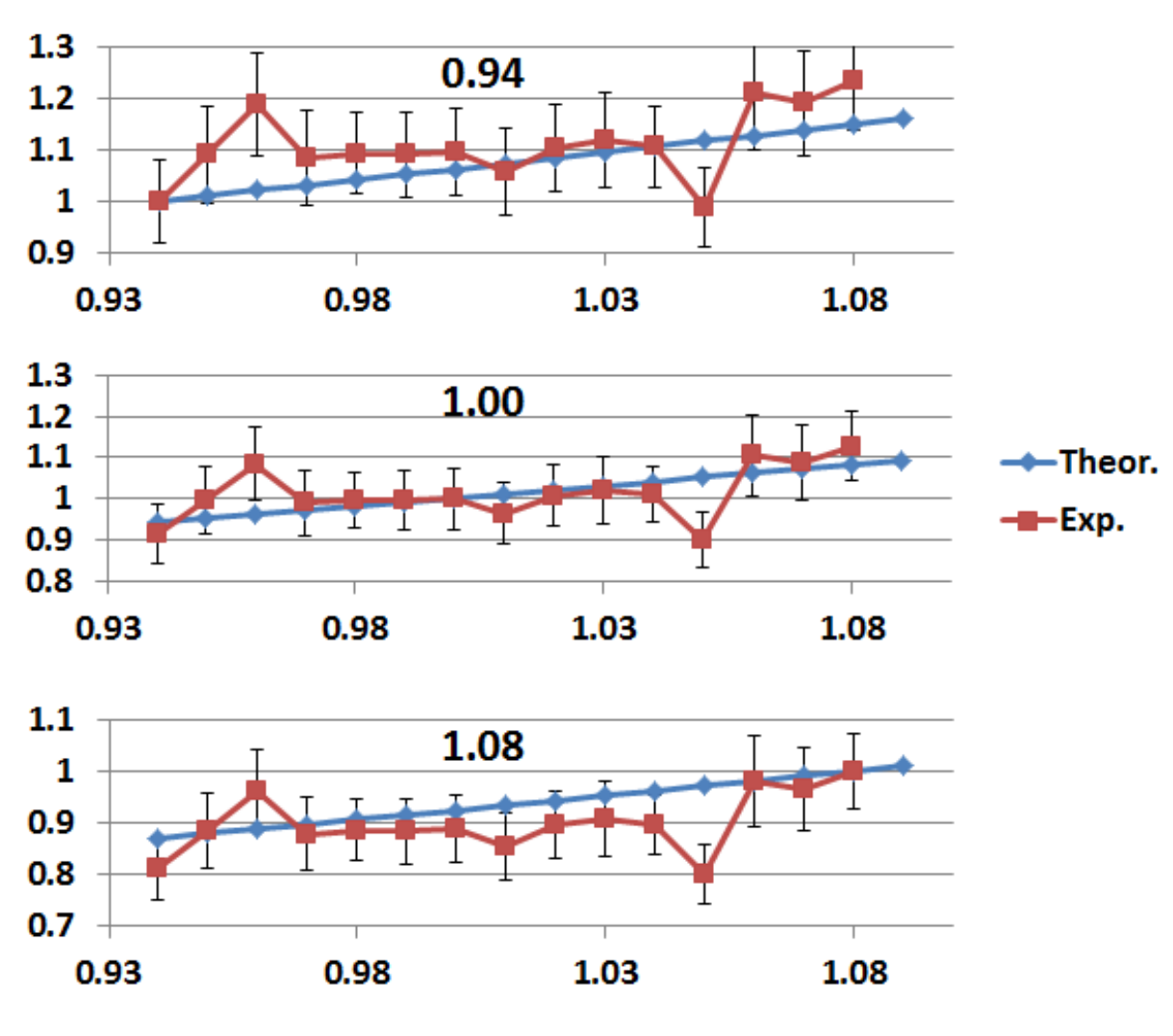}
\caption{The comparison between the theoretical $R_m$ ratios (computed from Equation~\ref{RatioValanga}) and the experimental ones, the subset of results corresponding to $m=1$ and maintaining fixed $\lambda_a$ respectively to 0.94, 1.00 and 1.08. The bars correspond to the standard deviation associated to each experimental $R_m$ quotient: they are computed by propagating through the quotient the standard deviations associated to each experimental avalanche distribution, estimated by means of bootstrapping procedures. The experimental bins have a size of $0.01$, each bin grouping from 110 till a maximum of 997 knock-out events.}
 \label{ratios}
\end{figure}

\begin{figure} 
\centering
 \includegraphics[width=0.8\textwidth]{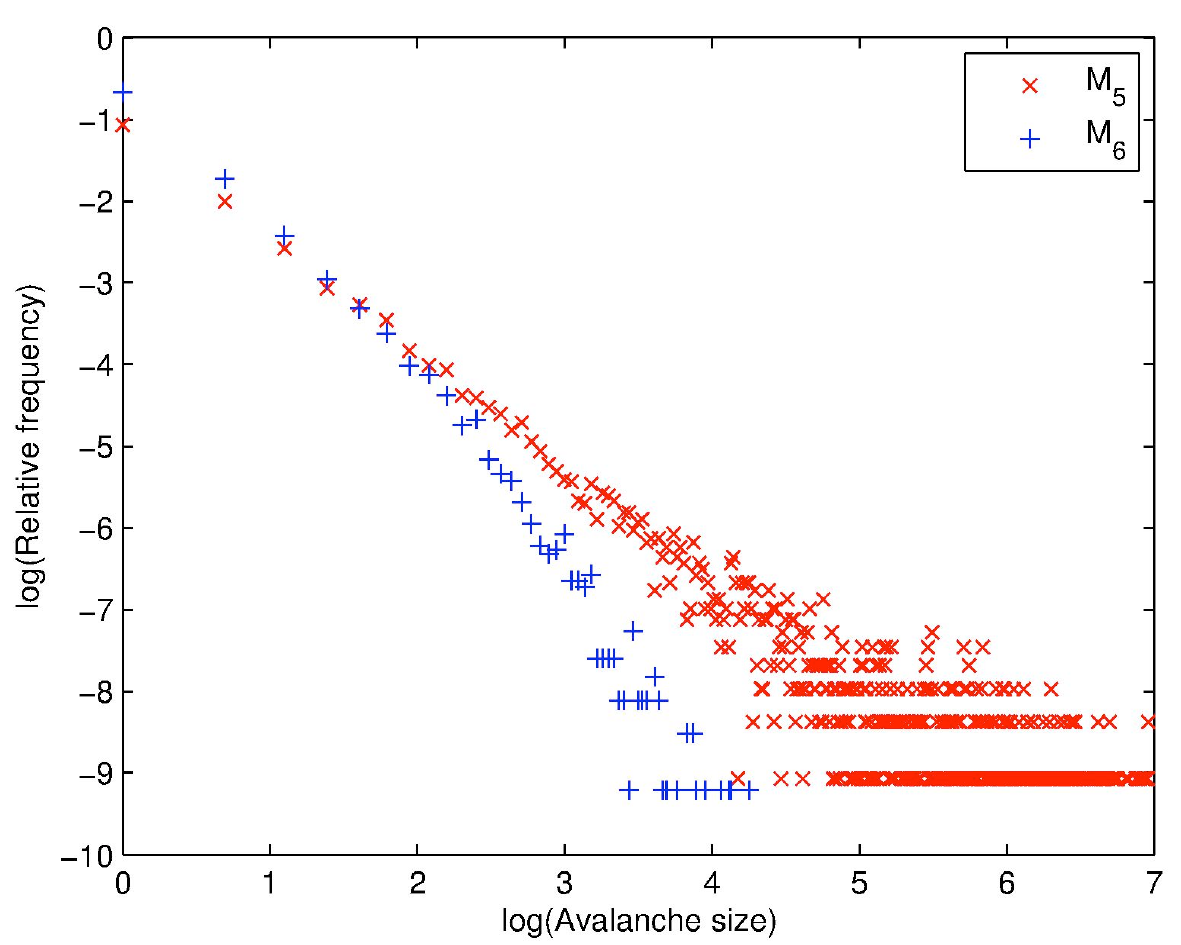}
\caption{Avalanche distribution in $\mathcal{M}_5$  (x) and $\mathcal{M}_6$ (+) families. It is possible to observe the difference on the distribution tails. $N=6000$, 10000 networks for each family.}
 \label{avalanchesLog}
\end{figure}

In~\cite{serra2004} it was shown that  RBN families with random topology, $N=6300$ (the number of genes in {\em Saccharomyces cerevisiae}), $k_{in}=2$ and a particular set of Boolean functions (all the 13 two-input Boolean functions that remain after excluding the not biologically plausible $XOR$, $EQUIVALENCE$ and $FALSE$ functions). These choices allowed a very close approximation of the experimental distributions. The measures we are proposing in this article, when applied to these networks, give respectively $DA=0.91$ and $SA=0.94$. These values of $DA$ and $SA$ are not very different, so if we were to use the former in Equation~\ref{valanga} we would obtain nonetheless a good approximation of the avalanche size distribution.
In our previous papers on avalanche distribution~\cite{serra2004,serra2007}, we were interested in comparison with available experimental data. In this case no information is available about the features of the attractor state of the yeast, therefore it is correct to use average values such as $DA$ and $SA$. However, in synthetic networks we can test the effect of the different attractor sensitivities on the resulting avalanches.
By applying Equation~\ref{valanga} to avalanches having different $SA_i$, (i.e. different $\lambda$'s) we obtain that the ratio between the frequency of avalanches of the same size $m$, starting from attractor $a$ or attractor $b$, is ruled by the following equation:
\begin{equation}
\label{RatioValanga}
R_{m}^{a,b}=\frac{ B_m\lambda_a^{m-1}e^{-m\lambda_a}  }{B_m\lambda_b^{m-1}e^{-m\lambda_b}}=(\frac{\lambda_a}{\lambda_b})^{m-1} e^{-m(\lambda_a - \lambda_b)}
\end{equation}
Note that these ratios are independent from the coefficients $B_m$.\\
We checked this relationship by building 8000 RBN realizations, reaching one attractor for each realization (starting from random initial conditions), measuring its $SA_i$ value and performing a knock-out event. 
We also grouped the $SA_i$ values in bins, small but however sufficiently ample to have enough points in each bin. The comparison with the experimental results is good, as it can be seen in the three cases shown in Figure~\ref{ratios} (corresponding to the fixed $\lambda_a$ values 0.94, 1.00 and 1.08).
As shown in previous paragraphs, networks with a limited set of Boolean functions or evolved networks can have very different $DA$, $SA$ and $SA_i$: in these cases the effect on the avalanche distributions are even more evident. A clear example is provided by $\mathcal{M}_5$ and $\mathcal{M}_6$ families: in these ensembles the theoretical and experimental $DA$ are respectively equal to 0.75 and 0.74, but the SA measures are equal to 0.96 and 0.67, these values originating from asymptotic states having respectively $b=0.08$ and $b=0.67$. These different asymptotic dynamical regimes give rise to very different avalanche distributions, as shown in figure~Figure~\ref{avalanchesLog}  (note the different distribution tails). In this case the classical Derrida parameter $DA$ would provide no information concerning the avalanche distribution.

\newpage
\section{Conclusion}
\label{Conclusions}
RBNs are widely used to simulate genetic regulatory systems, and it is well-known that they can show different dynamical regimes. In this paper we introduced different measures related to these dynamical behaviours, and we have shown that they provide a richer description of the dynamics than the usual ones, based on static sensitivity or on perturbations on random initial states. While some of the remarks presented here can be found also in other papers, which we have reviewed in this contribution, this work provides a more comprehensive analysis, and shows examples of a number of network families with somewhat unexpected properties.

The main claims of this paper can be summarized as follows:
\begin{itemize}
  \item Like other dynamical systems, RBNs can show ordered, disordered and critical regimes. The same RBN can show several different asymptotic behaviours (attractors);
  \item three kinds of measure are useful: {\em (i)} the classical Derrida measure ($DA$ in the text), which characterizes each RBN structural parameters, {\em (ii)} the attractor sensitivity ($SA_i$ in the text, relative to attractor $i$), which characterizes the dynamical regime of each attractor, and {\em (iii)} the attractors sensitivity ($SA$ in the text), which is related to the whole set of asymptotic dynamical behaviours of a RBN;
  \item in classically critical RBNs $DA$ and $SA$ are very close, and they tend to coincide when considering ensemble averages;
  \begin{itemize}
    \item single RBN realizations can show asymptotic behaviours with attractor sensitivity considerably far from its $DA$ value;
    \item and this tendency grows as the node connectivity $k_{in}$ increases;
  \end{itemize}
  \item finally, by knowing the asymptotic proportion of 0s and 1s it is possible to relate static and dynamic measures.
\end{itemize}

We remark that the issues raised in this paper are not limited to RBNs, but they hold for a large class of discrete computational models and they can have strong consequences also for the biological processes that they describe. Indeed, these systems spend most of their time in their attractor cycles, so the most informative stability analysis should be performed in those regimes. The most common approach, based on perturbing random initial states, might fail in properly describing the asymptotic behaviour of these non-ergodic systems, except for some particular conditions.

\section*{Authors' contributions}

MV and RS conceived and designed the experiments and provided a first analysis of the results. DC developed the code and performed the experiments. CD performed further experiments. MV, RS, CD, AR and AF analysed and discussed the results. AR, MV, CD, and AF wrote the paper. All authors gave final approval for publication.

\section*{Acknowledgements}
We are deeply indebted to Stuart Kauffman for his inspiring ideas and for several discussions on various aspects of random Boolean networks. We also gratefully acknowledge useful discussions with David Lane and Alex Graudenzi.

\bibliographystyle{plain}
\bibliography{biblio-campioli}

\end{document}